\documentstyle[12pt]{article}

  \newcommand{\const}{\rm const}

\begin{document}

   \begin{center}

{\bf Grand Lebesgue Spaces norm estimates for eigen functions  for }\\

\vspace{4mm}

{\bf Laplace - Beltrami operator defined on the }

\vspace{4mm}

{\bf closed compact smooth Riemannian manifolds.} \par

 \vspace{5mm}

 {\bf M.R.Formica,  E.Ostrovsky and L.Sirota. } \par

\vspace{5mm}

\end{center}

 Universit\`{a} degli Studi di Napoli Parthenope, via Generale Parisi 13, Palazzo Pacanowsky, 80132,
Napoli, Italy. \\

e-mail: mara.formica@uniparthenope.it \\

 \vspace{4mm}

Department of Mathematics and Statistics, Bar-Ilan University, \\
59200, Ramat Gan, Israel. \\

\vspace{3mm}

e-mail:eugostrovsky@list.ru\\
Department of Mathematics and Statistics, Bar-Ilan University,\\
59200, Ramat Gan, Israel.

\vspace{3mm}

e-mail:sirota3@bezeqint.net \\

\begin{center}

  {\bf Abstract} \par

\vspace{4mm}

 \end{center}

 \hspace{3mm}  We derive a sharp Grand Lebesgue Space norm estimations for normalized eigen functions
 for the Laplace - Beltrami operator defined on the compact smooth  Riemann manifold. \par
  \ These estimates allow us to deduce in particular the exponential decreasing tail of distribution for 
 these eigen functions. \par

\vspace{4mm}

\begin{center}

\vspace{4mm}

 \ {\it Key words and phrases:}

\end{center}

  \hspace{3mm}

 \ Compact smooth closed Riemann manifold,  Laplace - Beltrami operator, eigen values and functions, Lebesgue - Riesz
 and Grand Lebesgue norms and spaces,  tail function, Young - Fenchel transform, Young inequality, fundamental function,
 subgaussian variables, normalized function, non - asymptotic upper and lower estimate, generating function.

 \vspace{5mm}

 \section{Statement of problem. Notations. Previous results.}

 \vspace{5mm}

  \hspace{3mm} Let $ \ (M, g) \ $ be a compact closed smooth Riemannian manifold of dimension  $ \ d \ge 2, \ $ and let
  $ \ \Delta_g \ $ be the associated Laplace - Beltrami operator.   We will consider the $ \ L_2 \ - \ $ normalized eigenfunctions
satisfying  the classical relations

\begin{equation} \label{eigen equations}
- \Delta_g e_{\lambda}(x) = \lambda^2 e_{\lambda}(x), \ ||e_{\lambda}||_2^2 = \int_M |e_{\lambda}(x)|^2\ V_g(dx) = 1, \ \lambda > 0,
\end{equation}
where $ \ V_g(dx) \ $  (measure)  is element of volume on $ \ M \ $  and as ordinary $ \ ||f||_p \ $ denotes the classical Lebesgue - Riesz norm
for the (measurable) function $ \ f: M \to R: \ $

$$
||f||_p \stackrel{def}{=} \left[ \ \int_M |f(x)|^p \ V_g(dx) \ \right]^{1/p}, \ p \ge 2.
$$

 \ Introduce the following variables

$$
p_c := \frac{2(d + 1)}{d-1}, \ d \ge 2;
$$

$$
\mu(p) := \frac{d-1}{2} \cdot \left( \ \frac{1}{2}   - \frac{1}{p}  \ \right), \hspace{2mm} 2 < p \le p_c;
$$

$$
\mu(p) := d \left( \frac{1}{2} - \frac{1}{p} \right) - \frac{1}{2}, \hspace{2mm} p_c \le p \le \infty.
$$

 \ We will apply the following important estimate

\begin{equation} \label{source relation}
||e_{\lambda}||_p \le C(M,g) \ \lambda^{\mu(p)}, \ p > 2,
\end{equation}
see \cite{Sogge 1}, \cite{Sogge 2} and another works of this author \cite{Sogge 3} - \cite{Sogge 9}.
See also the articles \cite{Donnelly},   \cite{Xu 1} - \cite{Xu 2}.\par

\vspace{5mm}

 \hspace{3mm} {\bf  We intent to extend  the estimate (\ref{eigen equations}) from the classical Lebesgue - Riesz spaces
 into the more general  ones, namely, into the so - called Grand Lebesgue Spaces. } \par

\vspace{5mm}

\begin{center}

{\sc A brief review of the theory of  Grand Lebesgue Spaces.}

\end{center}

\vspace{5mm}

 \hspace{3mm} Let $ \ (a,b) = \const, \ 1 \le a < b \le \infty, \ $ and let $ \ \psi = \psi(p), \ p \in (a,b) \ $  be
 bounded from below: $ \ \inf_{p \in (a,b) } \psi(p) > 0 \ $ measurable function. The set of all such a functions will be
 denoted by  $ \ \Psi(a,b); \ $ put also

$$
\Psi  := \cup_{(a,b): 1 < a  < b < \infty} \Psi(a,b).
$$

\vspace{4mm}

 \ {\bf Definition 1.1.}  Recall that the so - called
 Grand Lebesgue Space $ \ G\psi, \ \psi \in \Psi(a,b) \ $  builded
 in particular on the set $ \ M \ $ equipped as before with the measure $ \ V_g, \ $
  consists by definition on all the integrable  numerical valued functions having  a finite norm

\begin{equation} \label{GLS norm Rd}
||f||G\psi =  ||f||G\psi(M) \stackrel{def}{=} \sup_{p \in (a,b)} \left\{ \ \frac{||f||_p}{\psi(p)}   \ \right\}.
\end{equation}

\vspace{4mm}

 \ The function $ \  \psi = \psi(p), \ p \in (a,b)  \ $ is said to be as {\it generating function} for this space. \par

\vspace{3mm}

 \ These spaces are rearrangement invariant Banach functional spaces. They
  was investigated in many works, see e.g. \cite{Buldygin}, \cite{Ermakov etc. 1986},
\cite{Fiorenza 1}, \cite{Fiorenza 2}, \cite{Fiorenza-Formica-Gogatishvili-DEA2018},
\cite{fioforgogakoparakoNAtoappear}, \cite{fioformicarakodie2017}, \cite{formicagiovamjom2015}, \cite{Iwaniec},
\cite{Kozachenko 1}, \cite{Kozachenko 2}, \cite{Ostrovsky 0} - \cite{Ostrovsky 3}. In particular, the belonging of
the function to certain Grand Lebesgue Space $ \ G\psi \ $ is closely related with its tail behavior and  is related
with its moment generating function

$$
\nu[f](z) := \int_M \exp( z \ f(x)) \ V_g(dx).
$$

 \ Notice that if we choose the following {\it extremal}  function

$$
 \psi_r = \psi_r(p)= 1, \ r = p, \  \psi_r(p) = \infty, \ p \ne r, \  \ r = \const > 1,
$$
and agree to take $ \ C / \infty = 0, \ $ then

$$
||f||G\psi_r = ||f||_r.
$$
 \ So, the notion of GLS contains as a particular case the classical Lebesgue - Riesz one. \par
 \ Further, let $ \ f(\cdot) \in G\psi\ $ and (for definiteness) $ ||f||G\psi = 1. \ $ Define  the following function
 (Young - Fenchel transform)

$$
h[\psi](u) := \sup_{p \in (a,b)} (pu - p \ln \psi(p)), \ u \ge e.
$$
 \ Then the tail function $ \ T[f](u) \ $ for $ \ f(\cdot) \ $ may be estimated as follows

$$
T[f](u) \stackrel{def}{=} V_g \{x, \ x \in M, \ |f(x)| > u \ \} \le \exp( - h[\psi](u)), \ u \ge e,
$$
the exponential decreasing in general case estimate; and inverse conclusion holds true up to finite constant under appropriate
natural conditions.\par

 \vspace{4mm}

 \ The fundamental function for these spaces  $ \ \phi[G\psi](\delta) = \phi[G\psi(a,b)](\delta), \ \delta > 0  \ $ has a form

\begin{equation} \label{fund fun}
\phi[G\psi(a,b)](\delta) = \sup_{p \in (a,b)} \left\{ \ \frac{\delta^{1/p}}{\psi(p)} \ \right\}.
\end{equation}
 \ These function was investigated in particular in \cite{Ostrovsky 3}. They used in functional analysis, theory of Fourier series etc. 
 They are also closely and continuously related with generating function $ \ \psi(p). \ $ \par
  \ A very important particular subgaussian case: $ \ \psi(p) = \sqrt{p}, \ p \in (1,\infty). \ $  \par

 \vspace{5mm}

 \section{Main result.}

 \vspace{5mm}

 \hspace{3mm} {\bf Case A: \ small values of the parameter.} \par

 \vspace{3mm}

  \ We consider here at first the case when $ \ 2 < p  \le p_c. \ $ Let $ \ 2 < a < b \le p_c \ $  and let $ \ \psi \in \Psi(a,b). \ $  \par

\vspace{5mm}

\ {\bf Theorem  2.1.}

\vspace{3mm}

\begin{equation} \label{small p}
||e_{\lambda}||G\psi \le C(M,g) \ \lambda^{(d - 1)/4} \ \phi[G\psi] \left(\lambda^{(1 - d)/2} \right), \ \lambda > 0.
\end{equation}

\vspace{5mm}

\ {\bf Proof.} We have from the source relation (\ref{source relation}) taking into account the restrictions $ \ \lambda > 0, \ 2 < p \le p_c \ $

$$
||e_{\lambda}||_p \le C(M,g) \ \lambda^{(d-1)/4} \ \lambda^{(1 - d)/(2p)},
$$
and after dividing over $ \ \psi(p) \ $

$$
\frac{||e_{\lambda}||_p}{\psi(p)} \le C(M,g) \ \lambda^{(d-1)/4} \cdot \frac{(\lambda^{(1 - d)/2})^{1/p}}{\psi(p)}.
$$
\ It remains to take the supremum over $ \ p, \ p \in (a,b) \ $ to get (\ref{small p}).\par

\vspace{4mm}

\hspace{3mm} {\bf Case B: \ great values of the parameter.} \par

 \vspace{3mm}

  \ Let us consider now  the case when $ \  p \ge  p_c. \ $ Let  here $ \ p_c \le a < b \le \infty \ $  and let $ \ \psi \in \Psi(a,b). \ $  \par

\vspace{5mm}

\ {\bf Theorem  2.2.}

\vspace{3mm}

\begin{equation} \label{great p}
||e_{\lambda}||G\psi \le C(M,g) \ \lambda^{(d - 1)/2} \ \phi[G\psi] \left(\lambda^{ - d} \right), \ \lambda > 0.
\end{equation}

\vspace{5mm}

 \ {\bf Proof } is quite alike ones in the foregoing case. Namely, we have for the values $ \ p \ge p_c \ $

$$
\mu(p) = \frac{d-1}{2} - \frac{d}{p},
$$
following in this case

$$
||e_{\lambda}||_p \le C(M,g) \ \lambda^{(d-1)/2} \ \lambda^{ - d/p},
$$

$$
\frac{||e_{\lambda}||_p}{\psi(p)} \le C(M,g) \cdot \lambda^{(d-1)/2} \times  \frac{(\lambda^{ - d})^{1/p}}{\psi(p)},
$$
which follows in turn to (\ref{great p}) after taking the supremum over $\ p. \ $ \par

\vspace{4mm}

 \ {\bf Example 2.1.} Let us choose $  \ \psi(p) = 1, \ p > p_c; \ $ and one can take $ \ p\to \infty; \ $ then we conclude

$$
\max_{x \in M} |e_{\lambda}(x)| \le C(M,g) \cdot \lambda^{(d-1)/2}, \ \lambda > 0.
$$

 \vspace{6mm}

\vspace{0.5cm} \emph{Acknowledgement.} {\footnotesize The first
author has been partially supported by the Gruppo Nazionale per
l'Analisi Matematica, la Probabilit\`a e le loro Applicazioni
(GNAMPA) of the Istituto Nazionale di Alta Matematica (INdAM) and by
Universit\`a degli Studi di Napoli Parthenope through the project
\lq\lq sostegno alla Ricerca individuale\rq\rq .\par

\end{document}